\documentclass[a4paper]{article}

\usepackage{color}
\usepackage{yfonts}
\usepackage{graphicx}
\usepackage{amssymb}
\usepackage{amsmath}
\usepackage[all]{xy}
\usepackage{amscd}
\usepackage{MnSymbol}
\usepackage{calc}
\usepackage{hyperref}

\newcommand{\mbr}{\mathbb{R}}
\newcommand{\mbd}{\mathbb{D}}

\newcommand{\mbsr}{\hat{\mathbb{C}}}
\newcommand{\mbn}{\mathbb{N}}
\newcommand{\mbz}{\mathbb{Z}}
\newcommand{\mbq}{\mathbb{Q}}

\newtheorem{theorem}{Theorem}[section]

\newtheorem{lemma}[theorem]{Lemma}

\newtheorem{corollary}[theorem]{Corollary}

\newtheorem{proposition}[theorem]{Proposition}

\newtheorem{property}[theorem]{Property}

\newtheorem{properties}[theorem]{Properties}

\newtheorem{example}[theorem]{Example}

\newtheorem{definition}[theorem]{Definition}

\newcommand{\pr}[1]{\vspace{1mm}\noindent\textit{Proof #1}.}
\newcommand{\finpr}{$\blacksquare$ \vspace{3mm}}
\newcommand{\finsspr}{$\square$ \vspace{3mm}}

\newcommand{\vs}[1]{\vspace{#1mm}}

\newcommand{\ju}[0]{\mathcal{J}}

\newcommand{\inter}{\partial B_1\cap\partial B_2}

\newcommand{\oc}[0]{\mathcal{P}}

\newcommand{\chis}{\chi(S^2,GO(\oc))}

\title{Periodic points in the intersection of immediate attracting basins boundaries}
\date{}
\author{Bastien Rossetti}

\begin{document}

\maketitle

\begin{center}
\begin{minipage}{13cm}
\footnotesize{\textbf{Abstract :} We give sufficient conditions under which the set of eventually periodic points in the intersection of immediate attracting basins boundaries is non-empty. We give other conditions under which this set is dense in the intersection.}
\end{minipage}
\end{center}

\vs{10}
For a rational map $f$ of degree $D\geq2$ with two distinct attracting periodic points $a_1$ and $a_2$, we wonder if $\inter$ contains a periodic point, where $B_i$ is the connected component of the complementary of the Julia set $\ju$ of $f$ containing $a_i$. In this note we give conditions to a positive answer.

To state the theorem, let $\oc:=\{f^n(c) : f'(c)=0~; n\geq1\}$ be the postcritical set of $f$ and let $\oc_b$ be the set of $x\in\oc$ such that there is no component of $\mbsr-\ju$ whose boundary contains $x$.

\begin{center}
\begin{minipage}{13cm}
\textbf{Theorem} Suppose that the $B_i$ are simply connected, the boundaries $\partial B_i$ are locally connected and that $\inter\neq\emptyset$.
\begin{enumerate}
\item[1.] If $x\in\inter$ is multiply accessible from $B_1$ and if the orbit of $x$ does not contain any critical point, then $x$ is eventually periodic.
\item[2.] If $\inter$ contains no critical point with infinite orbit and is disjoint from the $\omega$-limit set of every recurrent critical point, then $\inter$ contains a periodic point.
\item[3.] Assume furthermore that $\ju$ is connected. If $\inter$ contains no accumulation point of $\oc_b$ and $\oc\cap{(\partial B_1\cup\partial B_2)}$, then the subset of eventually periodic points in $\inter$ is non-empty and dense in $\inter$.
\end{enumerate}
\end{minipage}
\end{center}

As a particular case if $f$ is postcritically finite then the set of eventually periodic points in $\inter$ is non-empty (and dense in $\inter$). Nevertheless the theorem does not require $\oc$ to be finite.

The hypothesis of part 2 are not useless for the existence of a periodic point in $\inter$. Pick $\theta\in\mbr/\mbq$ and consider $F_\theta(z)=\rho_\theta z^2(z-3)/(1-3z)$ with $\rho_\theta\in S^1$ such that $F_\theta$ has rotation number $\theta$ on $S^1$. Such a map is studied in \cite{Pet}. The map $F_\theta$ has two attracting points $a_1=0$ and $a_2=\infty$, with $1\in\inter\subset S^1$. Altought the $B_i$ are simply connected and the $\partial B_i$ are locally connected there is no periodic point in $\inter$, since $F_{\theta|S^1}$ is topologically conjugate to $z\mapsto e^{2i\pi\theta}z$. Note that for these examples the point $1$ is a critical point in $\inter$ with an infinite orbit.

\vs{2}
For the proof of the theorem we assume that the $B_i$ are fixed, up to consider an iterate of $f$, and that $a_i$ is the only critical point of $f$ in $B_i$, up to make a quasi-conformal deformation (see \cite{CG}, theorem VI 5.1). Let $\phi_i:\mbd\to B_i$ conjugate $f$ to $z^{d_i}$ (where $d_i$ is the local degree of $f$ at $a_i$), and for $t\in\mbr$ let $R_i(t):=\phi_i(\{re^{2i\pi t} : 0\leq r<1\})$. Since $\partial B_i$ is locally connected, $\phi_i$ extends continuously to $\overline{\phi_i}:\overline{\mbd}\to\overline{B_i}$. A point $y\in\partial B_i$ is \textit{multiply accessible from $B_i$} if $\sharp((\overline\phi_i)^{-1}(\{y\}))\geq2$.

\vs{2}
\pr{of part 1}
Let $\theta,\theta'\in(\overline{\phi_1})^{-1}(\{x\})$ be distinct. For any $n\geq0$, $R_1(d_1^n\theta)$ and $R_1(d_1^n\theta')$ are distinct, for otherwise the orbit of $x$ contains a critical point. The union $\overline{R_1(d_1^n\theta)}\cup\overline{R_1(d_1^n\theta')}$ is a Jordan curve whose complementary has one connected component $V_n$ such that $\overline{V_n}\cap\overline{B_2}=\{f^n(x)\}$. Let $I_n:=(\overline{\phi_1})^{-1}(\partial B_1\cap\overline{V_n})$. If there exist $M\neq N$ with $\{d_1^N\theta,d_1^N\theta'\}\cap I_M\neq\emptyset$, then $f^M(x)=f^N(x)$. Otherwise the $I_k$ are disjoint intervals. This is impossible since $t\mapsto t^{d_1}$ is expanding on the unit circle, so that there are infinitely many $I_n$ with $|I_n|\geq1/d^2_1$.
\finpr

For the proof of part 2 we show (under hypothesis) that $f_{|\inter}$ is distance-expanding using a theorem of Ma\~{n}\'e and we conclude with a result about periodic points for distance-expanding maps.

\begin{theorem}\label{mane}
\textnormal{\textbf{(\cite{Ma})}} Let $g$ be a rational map of degree at least $2$, and $\Lambda$ be a compact forward invariant set in the Julia set containing no critical point nor parabolic point. If $\Lambda$ is disjoint from the $\omega$-limit set of every recurrent critical point, then there exists $N\geq0$ such that $\min_{z\in \Lambda}||(g^n)'(z)||>1$ for every $n\geq N$.
\end{theorem}

\noindent Let $(X,\rho)$ be a metric space. A continuous map $T:X\to X$ is \textit{distance-expanding} if there exist $\lambda>1$, $\eta>0$ and $N\geq0$ such that for any $x,y\in X$, if $\rho(x,y)\leq\eta$ then $\rho(T^N(x),T^N(y))\geq\lambda \rho(x,y)$.

\begin{theorem}\label{prur}
\textnormal{\textbf{(\cite{PrUr}, chapter 4)}} Let $(X,\rho)$ be a compact metric space. If $T:X\to X$ is continuous, open and distance-expanding, then there exists $\alpha>0$ such that if $x\in X$ and $\rho(x,T^L(x))\leq\alpha$ for some $L\geq1$, then $X$ contains a periodic point.
\end{theorem}

\pr{of part 2} Assume that $\inter$ contains no critical point with finite orbit nor parabolic point, nor multiple point in $\partial B_i$ (for otherwise $\inter$ contains a periodic point). By theorem~\ref{mane} there exists $N\geq0$ such that $\min_{x\in\inter}||(f^N)'(x)||>1$. We show now that $f_{|\inter}$ is distance-expanding and open, in order to apply theorem \ref{prur}.

\begin{lemma}\label{definitionsequivalence}
The restriction $f_{|\inter}$ is distance-expanding with respect to the spherical metric.
\end{lemma}

\pr{}
By continuity of $x\mapsto||(f^N)'(x)||$, there exist $\lambda>1$ and a neighborhood $U$ of $\inter$ such that $\min_{x\in U}||(f^N)'(x)||\geq\lambda$. By compactness of $\inter$ there exists $\eta>0$ such that if $d(x,y)\leq\eta$ then the geodesic $\Gamma$ between $f^N(x)$ and $f^N(y)$ lifts to a path $\gamma$ from $x$ to $y$ with $\gamma\subset U$. Thus we get $d(f^N(x),f^N(y))=\textnormal{length}(\Gamma)\geq\lambda.\textnormal{length}(\gamma)\geq\lambda d(x,y)$.
\finsspr

\begin{lemma}
The restriction $f_{|\inter}$ is open.
\end{lemma}

\pr{}
Let $O\subset\inter$ and suppose $f(O)$ is not open. There exists a sequence $(y_n)$ in $(\inter)-f(O)$ tending to some $y\in f(O)$. Let $x\in O$ be such that $f(x)=y$. There exist a neighborhood $U$ of $x$ and a neighborhood $V$ of $y$ such that $f:U\to V$ is a homeomorphism. Thus for $n$ large enough the point $x_n=f^{-1}(y_n)\cap U$ is well defined, and $x_n\to x$. We show now that $x_n\in \inter-O$ and get that $O$ is not open. It is clear that $x_n\notin O$ since $f(x_n)\notin f(O)$. For any $n$ there exists a Fatou component $B_i^n$ such that $f(B_i^n)=B_i$ and $x_n\in\partial B_1^n\cap\partial B_2^n$. For $n$ large enough $B_i^n=B_i$. Otherwise, for some $i_0\in\{1,2\}$ there exists a Fatou component $B$ such that $B\neq B_{i_0}$, $f(B)=B_{i_0}$ and $x\in\partial B$. The boundary $\partial B$ has finitely many connected components, thus each one of them is locally connected. Let $\tilde B$ be either $B$ or $B_{i_0}$. There exists a connected component $U_{\tilde B}$ of $U\cap\tilde B$ such that $x\in\partial U_{\tilde B}$. Since $f(x)$ is simple in $\partial B_{i_0}$, there exists a unique connected component $V_{B_{i_0}}$ of $V\cap B_{i_0}$ such that $f(x)\in\partial V_{B_{i_0}}$. Hence $f(U_{\tilde B})=V_{B_{i_0}}$. Since $B\neq B_{i_0}$, we have $U_{B_{i_0}}\cap U_B=\emptyset$ and $f(U_{B_{i_0}})=f(U_B)$ which contradicts the injectivity of $f_{|U}$.
\finsspr

\noindent Now we apply theorem \ref{prur}. Let $w$ be an accumulation point of the orbit of some point $z$ in $\inter$. There exist $P>Q\geq0$ such that $f^P(z)$ and $f^Q(z)$ are both in $B(w;\alpha/2)$, hence $d(f^Q(z),f^P(z))=d(f^Q(z),f^{P-Q}(f^Q(z)))\leq\alpha$, where $\alpha$ is the constant in theorem \ref{prur}. Thus we get a periodic point in $\inter$. $\blacksquare$

The proof of part 3 uses ideas and techniques developed by Kevin Pilgrim in his thesis (\cite{Pil}). In case where $f$ is postcritically finite and hyperbolic, the part 3 is a corollary of his work.

\noindent Up to make a quasi-conformal deformation, we assume that all the critical points in the basin of $a_i$ have a finite orbit. We also assume that $\sharp\oc>2$ otherwise $f$ is conjugate to $z\mapsto z^{\pm D}$ and the conclusion follows.

For any $t,t'\in\mbr$, if $\overline{R_1(t)}\cap\overline{R_2(t')}\neq\emptyset$ then $\overline{R_1(t)}\cup\overline{R_2(t')}$ is called a \textit{chord}. Denote $\chi$ the set of chords. If $\alpha\in\chi$ is periodic then the point $\alpha\cap\ju\in\inter$ is periodic. We equip $\chi$ with the Hausdorff distance $d_H$, which turns it into a compact metric space. The construction of a periodic chord will use the theory of isotopies. For any chord $\alpha$ and any set $X\subset\mbsr$, $[\alpha]_{X}$ will denote the isotopy class of $\alpha$ rel $X$. For any $\alpha,\beta\in\chi$ distinct the complement $\mbsr-(\alpha\cup\beta)$ has at least two connected components and at most three (with points of $\ju$ in each one of them). Thus for any $m\geq0$, $[\alpha]_{f^{-m}(\oc)}=[\beta]_{f^{-m}(\oc)}$ if and only if one connected component of $\mbsr-(\alpha\cup\beta)$ contains all the points of $f^{-m}(\oc)-\{a_1,a_2\}$.
\\(*) Let $GO(\oc):=\{f^m(c) : f'(c)=0, m\in\mbz\}$. From the hypothesis the set $GO(\oc)\cap(\inter)$ is finite. We deduce that if a chord $\alpha$ satisfies $\alpha\cap\ju\in GO(\oc)$ then $\alpha\cap\ju$ is eventually periodic. We denote $\chis:=\{\alpha\in\chi : \alpha\cap\ju\notin GO(\oc)\}$. We study the dynamics of the chords in $\chis$.

The proof is as follows. For $\alpha\in\chis$ if the sequence $([f^j(\alpha)]_{f^{-1}(\oc)})_{j\geq0}$ is eventually cyclical then $\alpha$ is eventually periodic (lemma \ref{lemma5}). Otherwise, remarking that $([f^j(\alpha)]_\oc)_{j\geq0}$ contains twice the same element (lemma \ref{lemma3}) we will build a sequence $(\beta_n)\subset\chis$ by a series of adjustments (lemma \ref{lemma4}) such that $\beta_n$ converges (lemma \ref{lemma2}) to a chord $\beta$ with $\beta\cap\ju$ eventually periodic. The density part will follow from the fact that we can build $\beta$ as close as we want from $\alpha$.

Let $\alpha\in\chis$. A \textit{lift of $\alpha$} is the closure of a connected component of $f^{-1}(\alpha-\{a_1,a_2\})$. If a lift of $\alpha$ is a chord then it belongs to $\chis$. 

\begin{lemma}\label{lemma4}
Let $N\geq1$ and $\alpha,\beta_N\in\chis$ be such that $\beta_N\in[f^N(\alpha)]_{\oc}$. There exists a unique $\beta_0\in\chis\cap[\alpha]_{f^{-1}(\oc)}$, such that $f^N(\beta_0)=\beta_N$ and $[f^i(\beta_0)]_{f^{-1}(\oc)}=[f^i(\alpha)]_{f^{-1}(\oc)}$ for every $0\leq i\leq N-1$. More precisely, $\beta_0\in[\alpha]_{f^{-N}(\oc)}$.
\end{lemma}

\pr{}
Since $f^{N-1}(\alpha)$ is a lift of $\alpha$, there exists a unique lift $\beta_{N-1}$ of $\beta_N$ with $\beta_{N-1}\in[f^{N-1}(\alpha)]_{f^{-1}(\oc)}$. Moreover $\beta_{N-1}\in\chis$. In particular, $\beta_{N-1}\in[f^{N-1}(\alpha)]_{\oc}$. For each $1\leq i\leq N$ we construct inductively a unique $\beta_{N-i}\in[f^{N-i}(\alpha)]_{f^{-1}(\oc)}$ with $f^i(\beta_{N-i})=\beta_N$ and $f^k(\beta_{N-i})$ in $[f^k(f^{N-i}(\alpha))]_{f^{-1}(\oc)}$ for any $0\leq k\leq i-1$ (note that more precisely $\beta_{N-i}\in[f^{N-i}(\alpha)]_{f^{-i}(\oc)}$). Moreover $\beta_{N-i}\in\chis$.
\finpr

\begin{lemma}\label{lemma2}
$\forall\varepsilon>0$ $\exists N\geq0$ : $\forall\alpha,\beta\in\chi(S^2,GO(\oc))$, if $[\alpha]_{f^{-N}(\oc)}=[\beta]_{f^{-N}(\oc)}$ then $d_H(\alpha,\beta)<\varepsilon$. As a consequence, if $[\alpha]_{f^{-n}(\oc)}=[\beta]_{f^{-n}(\oc)}$ for every $n\geq0$ then $\alpha=\beta$.
\end{lemma}

\pr{}
Let $\varepsilon>0$. There exists a constant $\eta>0$ such that : for any $\alpha,\beta\in\chi(S^2,GO(\oc))$, if $d_H(\alpha,\beta)\geq\varepsilon$ then in two of the connected components of $\mbsr-(\alpha\cup\beta)$ lie an open ball centered at a point of $\ju$ and with radius $\eta$ (for otherwise we can construct a sequence $((\delta_n,\gamma_n))\subset\chis^2$ such that $d_H(\delta_n,\gamma_n)\geq\varepsilon$, and such that there are no two connected components of $\mbsr-(\delta_n\cup\gamma_n)$ in which lie an open ball centered at a point of $\ju$ with radius $\eta_n$, with $\eta_n\to 0$ ; taking a accumulation point $(\delta,\gamma)$ of $((\delta_n,\gamma_n))$, we have $d_H(\delta,\gamma)\geq\varepsilon$ and in two of the connected components of $\mbsr-(\delta\cup\gamma)$ lie an open ball centered at a point of $\ju$ and with fixed radius $\rho>0$ ; we see that for $N$ large enough in two of the connected components of $\mbsr-(\delta_N\cup\gamma_N)$ lie an open ball centered at a point of $\ju$ and with radius $\rho/2$ ; take $N$ large enough so that $\eta_N<\rho/2$, this is a contradiction). Since $\sharp\oc>2$, there exists $N\geq0$ (depending only on $\eta$) such that each one of these balls contains a point of $f^{-N}(\oc)$. Hence there is a point of $f^{-N}(\oc)$ in two connected components of $\mbsr-(\alpha\cup\beta)$, thus $[\alpha]_{f^{-N}(\oc)}\neq[\beta]_{f^{-N}(\oc)}$.
\finpr

\begin{lemma}\label{lemma5}
For $\alpha\in\chis$, if the sequence $([f^n(\alpha)]_{f^{-1}(\oc)})_{n=0}^\infty$ is cyclical then $\alpha$ is periodic.
\end{lemma}

\pr{}
Suppose that $[f^{n+Q}(\alpha)]_{f^{-1}(\oc)}=[f^{n}(\alpha)]_{f^{-1}(\oc)}$ for any $n\geq 0$. We apply lemma \ref{lemma4} to $f^{n+Q}(\alpha)\in[f^n(\alpha)]_{\oc}$ : there exists a chord $\beta_n\in[\alpha]_{f^{-1}(\oc)}$ such that $f^n(\beta_n)=f^{n+Q}(\alpha)$ and for all $0\leq i\leq n-1$, $[f^i(\beta_n)]_{f^{-1}(\oc)}=[f^{i}(\alpha)]_{f^{-1}(\oc)}$. This chord is $f^{Q}(\alpha)$. Thanks to lemma \ref{lemma4} we also have $f^Q(\alpha)=\beta_n\in[\alpha]_{f^{-n}(\oc)}$. Since this is true for any $n\geq0$, we conclude by lemma \ref{lemma2} that $f^{Q}(\alpha)=\alpha$.
\finpr

\begin{lemma}\label{lemma3}
For any $\alpha\in\chis$ there exist $M\neq N$ such that $[f^M(\alpha)]_\oc=[f^N(\alpha)]_\oc$.
\end{lemma}

\pr{}
Suppose that for any $m,n\geq0$ distinct we have $[f^m(\alpha)]_\oc\neq[f^n(\alpha)]_\oc$. Since $(\chi,d_H)$ is compact, there exists a chord $\beta\in\chi$ on which accumulate a sub-sequence $(f^{n(k)}(\alpha))_{k=0}^\infty$. Since for any $k,k'$ distinct at least two connected components of $\mbsr-(f^{n(k)}(\alpha)\cup f^{n(k')}(\alpha))$ contain a point of $\oc$, one can construct a non-stationary sequence $(z_n)\subset\oc$ which accumulates a point $z\in\beta$. Up to extraction, we have $(z_n)\subset\oc_b$ or $(z_n)\subset\oc\cap(\partial B_1\cup\partial B_2)$. Indeed, assume that $(z_n)\cap\oc-\oc_b$ is infinite. Up to extraction assume $(z_n)\subset\oc-\oc_b$. By construction, there is an infinite subset of $(z_n)$ whose elements are pairwise separated by chords. We deduce that $(z_n)\subset\oc\cap(\overline{B_1}\cup\overline{B_2})$ since there exist finitely many Fatou components $U_1,\dots,U_N$ such that $\oc-\oc_b\subset\overline{U_1}\cup\dots\cup\overline{U_N}$. Since we assume that $a_i$ is the only point of $\oc$ in $B_i$, we conclude that $(z_n)\subset\oc\cap(\partial B_1\cup\partial B_2)$. In particular, $(z_n)\subset\ju$ and $z=\lim_{n\to\infty}z_n=\beta\cap\ju\in\inter$. This contradicts the hypothesis of part 3.
\finpr

\pr{of part 3}
Let $\alpha$ be a chord. We have three cases.

\textit{Case 1 : $\alpha\cap\ju\in GO(\oc)$.} Thus $\alpha\cap\ju$ is eventually periodic, as explained in (*).

\textit{Case 2 : $\alpha\in\chi(S^2,GO(\oc))$ and $([f^n(\alpha)]_{f^{-1}(\oc)})_{n=0}^\infty$ is eventually cyclical.} Then $\alpha$ is eventually periodic by lemma \ref{lemma5}, hence the point $\alpha\cap\ju$ is eventually periodic.

\textit{Case 3 : $\alpha\in\chi(S^2,GO(\oc))$ and $([f^n(\alpha)]_{f^{-1}(\oc)})_{n=0}^\infty$ is not eventually cyclical.} We built from $\alpha$ a chord $\beta$ fitting case 1 or 2, as follows.
\\By lemma \ref{lemma3} there exist $N\geq0$ and $Q\geq1$ such that $[f^{N+Q}(\alpha)]_{\oc}=[f^N(\alpha)]_{\oc}$. Set $\beta_0:=f^N(\alpha)$. We apply lemma \ref{lemma4} to $\beta_0\in[f^Q(\beta_0)]_\oc$ : there exists a chord $\beta_1\in[\beta_0]_{f^{-Q}(\oc)}$ such that $$[f^i(\beta_1)]_{f^{-1}(\oc)}=[f^i(\beta_0)]_{f^{-1}(\oc)}$$ for any $0\leq i\leq Q-1$, and $f^{Q}(\beta_1)=\beta_0$. Thus $[f^{kQ+i}(\beta_1)]_{f^{-1}(\oc)}=[f^i(\beta_0)]_{f^{-1}(\oc)}$ for any $k\in\{0,1\}$ and any $0\leq i\leq Q-1$, and $\beta_1\in[f^{2Q}(\beta_1)]_{\oc}$.
\\We build inductively a sequence of chords $(\beta_q)_{q=0}^\infty$ such that :
\\(i) $\beta_{q+n}\in[\beta_q]_{f^{-2^qQ}(\oc)}$ for any $n\geq0$, and
\\(ii) $[f^{kQ+i}(\beta_q)]_{f^{-1}(\oc)}=[f^i(\beta_0)]_{f^{-1}(\oc)}$ for any $k\in\{0,\dots,2^q-1\}$ and $0\leq i\leq Q-1$.

\noindent The sequence $(\beta_q)_{q=0}^\infty$ converges to a chord $\beta$ whose point $\beta\cap\ju$ is eventually periodic. Indeed, the convergence follows from (i) and lemma \ref{lemma2}. The limit $\beta$ is a chord since $\chi$ is compact. If $\beta$ is not in $\chis$ then $\beta$ fits case 1. If $\beta\in\chi(S^2,GO(\oc))$, then we get for the limit $[f^{kQ+i}(\beta)]_{f^{-1}(\oc)}=[f^i(\beta_0)]_{f^{-1}(\oc)}$ for every $k\in\mbn$ and $0\leq i\leq Q-1$. Hence $([f^i(\beta)]_{f^{-1}(\oc)})_{i=0}^\infty$ is cyclical, and $\beta$ fits case~2.

We have obtained a chord $\beta\in[f^N(\alpha)]_{f^{-Q}(\oc)}$ such that the point $\beta\cap\ju$ is eventually periodic. We apply lemma \ref{lemma4} to $\beta\in[f^N(\alpha)]_{f^{-Q}(\oc)}$ : there exists a chord $\gamma\in[\alpha]_{f^{-(N+Q)}(\oc)}$ such that $\gamma\cap\ju$ is eventually periodic. Now we make $N$ tend to $\infty$ by applying lemma \ref{lemma3} to $f^m(\alpha)$ for $m$ as large as wished. Then we get a sequence of chords $(\gamma_n)\to\alpha$ with $\gamma_n\cap\ju$ eventually periodic. By local connectivity of the $\partial B_i$ we get that $\gamma_n\cap\ju$ converges to $\alpha\cap\ju$. This proves the density.
\finpr


\begin{thebibliography}{00}

\bibitem[CG]{CG} {\sc L. Carleson and T. W. Gamelin},
{\em Complex Dynamics.}
Second edition, Springer, 1995.

\bibitem[Ma] {Ma} {\sc R. Ma\~{n}\'e }, 
{\em On a theorem of Fatou}, Bol. Soc. Bras. Mat. vol. 24, p. 1-11 (1993)

\bibitem[Pet]{Pet} {\sc C.L. Petersen} ,
 {\em Local connectivity of some Julia sets containing a circle with an irrational rotation.}
Acta Math., 177 (1996), 163-224

\bibitem[Pil] {Pil} {\sc K. Pilgrim }, 
{\em Cylinders for iterated rational maps}, Thesis (1994)

\bibitem[PrUr] {PrUr} {\sc F. Przytycki and M. Urba\'nski }, 
{\em Conformal Fractals : Ergodic Theory Methods}, The London Mathematical Society, Lecture Note Series 371 (2010)

\end{thebibliography}
\end{document}